\documentclass[11pt,reqno,twoside]{amsart}
\usepackage{latexsym}
\usepackage{a4}
\usepackage{amsmath,amsthm}
\usepackage{amsfonts}
\usepackage[all]{xy}
\usepackage{lscape}
\usepackage{graphicx}
\usepackage{color}
\makeatletter
\def\LaTeX{\leavevmode L\raise.42ex
\hbox{\kern-.3em\size{\sf@size}{0pt}\selectfont A}\kern-.15em\TeX} \makeatother

\newtheorem{Def}{Definition}[section]
\newtheorem{cor}[Def]{Corolary}
\newtheorem{lem}[Def]{Lemma}
\newtheorem{obs}[Def]{Remark}
\newtheorem{prop}[Def]{Proposition}
\newtheorem{teo}[Def]{Theorem}
\newtheorem{ej}[Def]{Example}

\newcommand{\Hom}{\mbox{Hom}_{\mathcal{C}}(\tilde{T},}
\newcommand{\add}{\mbox{add}(\tau\tilde{T})}

\newcommand{\rd}{representation dimension }

\usepackage{color}

\begin{document}

\begin{abstract}
The aim of this work is to study the representation dimension of cluster tilted algebras.
We prove that the weak representation dimension of tame cluster tilted algebras is equal to three. We construct a generator module that reaches the weak representation dimension, unfortunately this module is not always  a cogenerator. We show for which algebras this module gives the representation dimension.
\end{abstract}
\vskip -.5cm

\title[On representation dimension of tame cluster tilted algebras]
{On representation dimension of tame cluster tilted algebras}

\author[Gonz\'alez Chaio]{Alfredo Gonz\'alez Chaio}
\address{Depto de Matem\'atica, Facultad de Ciencias Exactas y Naturales, Universidad Nacional de Mar del Plata,  7600, Mar del Plata, Argentina} \email{agonzalezchaio@gmail.com}

\author[Trepode]{Sonia Trepode}
\address{Departamento de Matem\'atica, Facultad de Ciencias Exactas y Naturales,
Funes 3350, Universidad Nacional de Mar del Plata, CONICET. 7600 Mar del
Plata, Argentina} \email{strepode@mdp.edu.ar}


\maketitle

\section{Introduction}

The representation dimension was introduced by Auslander \cite{Au} in the early seventies. Due to the
connection of arbitrary artin algebras with representation finite artin algebras, it was  expected that this notion would give a reasonable way of measuring how far an artin algebra is from being of representation
finite type. In fact, the representation dimension characterizes artin algebras of finite
representation type. Later, Iyama proved in \cite{I} that the representation
dimension of an artin algebra is always finite and Rouquier \cite{R} constructed examples of algebras with arbitrarily large representation dimension. On the other hand, the representation dimension is a measure, in some way, of the complexity of the morphisms of the module category. The purpose of our work is to study the representation dimension of tame cluster tilted algebras.
 \vskip .3cm
Cluster algebras were introduced by Fomin and Zelevinsky \cite{FZ}. Later,  Buan, Marsh, Reineke, Reiten and Todorov \cite{BMRRT} defined the cluster category and developed a tilting theory using a special class of
objects, namely the cluster tilting objects. In \cite{BMR2}, Buan, Marsh and Reiten
introduced  and studied the endomorphism algebras of such cluster-tilting objects in a cluster category. These algebras are known as cluster tilted algebras and are deeply connected to tilted  and hereditary algebras. They also prove that the module category of a cluster tilted algebra is a proper quotient of the cluster category.
This last fact motivates us to investigate the relationship between the module theory of cluster-tilted algebras and the module theory of hereditary algebras. In particular, one can obtain a lot of information about the module category of a cluster tilting algebra using the information provided by the torsion theory defined by the tilting module in the hereditary algebra.
\vskip .3cm
In \cite{GT} we have studied the representation dimension of cluster concealed algebras and we proved that its representation dimension is equal to three. The aim of this paper is to extend the techniques developed there to construct a generator of the module category of tame cluster-tilted algebras.
\vskip .3cm

Tame cluster tilted algebras are induced from tilting modules over tame hereditary algebras. Therefore, the torsion theory induced by a tilting module in module category of a tame hereditary algebra became a very useful tool. We study the representation dimension of tame cluster tilted algebras. More precisely, we prove that the weak representation dimension of tame cluster tilted algebras is equal to three. We construct a generator module which reaches the weak representation dimension, unfortunately the constructed module will not always be a generator. Finally, we show when this generator module gives the representation dimension of a tame cluster tilted algebra.

Recently, Garc\'{\i}a Elsener and Schiffler in \cite{GS} proved that the representation dimension of tame cluster tilted algebras is equal to three. The proof is based on a result of Bergh and Oppermann, that the category of Cohen Macaulay modules of this algebras is finite. Since the representation dimension is a measure of the complexity of the morphisms of the module category is still interesting to look for an Auslander generator of the module category. In this sense, the techniques introduced in this paper could be useful to know how the morphisms behave in the derived category.

In section 2, we introduce the notation we will use through the paper and recall some well known results about the Auslander-Reiten quiver of an hereditary algebra.
We also recall some preliminary result about torsion theory for hereditary algebras. In section 3, we
present our main theorem and prove various results required for
the proof of our main theorems.

\vskip .3cm

\section{preliminaries}

Through this paper, all algebras will be finite dimensional algebras over an algebraically closed
field and for an algebra $A$ we will consider mod $A$ the category of all finitely generated right $A$-modules and denote by ind$A$  the full subcategory of indecomposable modules in mod $A$.

Given an $A$-module $M$ we denoted by add$(M)$  the additive subcategory of mod $A$ generated by $M$, whose objects consist in finite sums of direct summands of the module $M$.
\vspace{0.2in}

We will denote by $H$ a finite dimensional hereditary
algebra and by $\mathcal{D}^{b}(H)$ the bounded derived category of $H$. Recall that the objects here are stalk complexes. We will identify the objects concentrated in degree zero with the corresponding $H$-modules. We will denote by [ ] the shift functor, thus any object in $\mathcal{D}^{b}(H)$ can be viewed as $X[i]$ with $X$ in mod $H$ and $i$ an integer.


\subsection{Representation dimension}
We recall that an $A$-module $M$ is a generator for mod $A$ if for
each $X \in$ mod $A$ there exists an epimorphism $M' \rightarrow X $
with $M' \in $ add$(M)$. Observe that $A$ is a generator for mod $A$.
Dually, we say that an $A$-module $M$ is a cogenerator if for each
$Y \in $ mod $A$ there exists a monomorphism $Y \rightarrow M'$ with
$M' \in $ add$(M)$. Note that $DA$ is a cogenerator for mod $A$. In
particular, any module $M$ containing every indecomposable
projective and e-very indecomposable injective module as a summand is
a generator-cogenerator module for mod $A$.

The original definition of representation dimension (we will denote it
by rep.dim) of an artin algebra $A$ is due to Auslander. For more
details on this topic, we refer the reader to \cite{Au}. The following
is a useful characterization of representation dimension, in the case
that $A$ is a non semisimple algebra, also due to Auslander.
The representation dimension of an artin algebra is given by

$$ \mbox{rep.dim}A = \mbox{inf} \,\{ \mbox{gl.dim End}_A(M) \mid M \mbox{is a generator-cogenerator for mod} A \}.$$

A module $M$ that reaches the minimum in the above definition is called an
$\mathbf{Auslander \, generator}$ and $\mbox{gl.dim End}_A(M)=
\mbox{rep.dim}A$ if $M$ is an Auslander generator.

The representation dimension can also be defined in a functorial
way, which will be more convenient for us. The next definition (see
\cite{APT},\cite{EHIS},\cite{CP},\cite{R}) will be very useful for the rest of
this work.

\begin{Def}\label{repdim}
The representation dimension $\mbox{rep.dim}$$A$ is the smallest integer $i
\geq 2$ such that there is a module $M \in \mbox{mod}A$ with the
property that, given any $A$-module $X$,

\begin{itemize}
\item [(a)] there is an exact sequence
$$0 \rightarrow M^{i-1} \rightarrow M^{i-2} \rightarrow ... \rightarrow M^{1}  \stackrel{f}\rightarrow X \rightarrow 0$$

\noindent with $M^{j} \in \mbox{add}(M)$ such that the sequence

$$ 0 \rightarrow \mbox{Hom}_{A}(M, M^{i-1}) \rightarrow ... \rightarrow  \mbox{Hom}_{A}(M, M^{1}) \rightarrow \mbox{Hom}_{A}(M, X) \rightarrow 0$$

\noindent is exact.
\item [(b)] there is a exact sequence
$$0 \rightarrow X \stackrel{g}\rightarrow M'_{1} \rightarrow M'_{2} \rightarrow ... \rightarrow M'_{i-1} \rightarrow 0$$

\noindent with $M'_{j} \in \mbox{add}(M)$ such that the sequence

$$ 0 \rightarrow \mbox{Hom}_{A}(M'_{i-1},M) \rightarrow ... \rightarrow  \mbox{Hom}_{A}(M'_{1}, M) \rightarrow \mbox{Hom}_{A}(X, M) \rightarrow 0 $$

\noindent is exact.
\end{itemize}
\end{Def}

Following \cite{CP}, we say that the module $M$ has the $i-1$-resolution property and that the sequence $0 \rightarrow M^{i-1} \rightarrow M^{i-2} \rightarrow ... \rightarrow M^{1}  \rightarrow X \rightarrow 0$ is an add$M$-approximation of $X$ of length $i-1$. Observe, that $f:M^1\to X$ is a right add$(M)$-approximation of $X$ and $g:X \rightarrow M'_1$ is a left add$(M)$-approximation of $X$.

\begin{obs}\label{1}
Either condition (a) or (b) imply that $\mbox{gl.dim End}_A(M) \leq i$ (\cite[Lemma 2.2]{EHIS}). Then, if $M \in \mbox{mod}A$ and $i \geq 2$, the following statements are equivalent.
\begin{itemize}
\item $M$ satisfies (a) and (b) of the definition.
\item $M$ satisfies (a) and $M$ contains an injective cogenerator as a direct summand.
\item $M$ satisfies (b) and $M$ contains a projective generator as a direct summand.
\end{itemize}
\end{obs}

\subsection{Torsion theory for hereditary algebras}
\vspace{0.3cm}

We start this section by giving the definition of tilting module over a hereditary algebra.
For more details on tilting modules see \cite{HR}. Let $H$ be an hereditary
algebra with $n$ non-isomorphic simple modules in mod$H$ and $T$ an $H$-module. $T$ is said to be a tilting module in
mod $H$ if it satisfies the following conditions:
\begin{itemize}
\item [(a)] Ext$^{1}_{A}(T, T) = 0$.
\item [(b)] $T$ has exactly $n$ non-isomorphic indecomposable direct summands.
\end{itemize}

\vskip .1cm

A tilting module is said to be basic if all its direct summands
are non-isomorphic. The endomorphism ring of a tilting module over a
hereditary algebra is said to be a $\mathbf{tilted}$ $\mathbf{
algebra}$. In particular, hereditary algebras are tilted algebras.
\vskip .1cm

We recall that a $\mathbf{path}$ from $X$ to $Y$ is a sequence $X = X_0 \rightarrow X_1 \rightarrow ...
\rightarrow X_t = Y$ with $t>0$ of non-zero non-isomorphisms between indecomposable modules. Given $X,Y \in \mbox{ind}A$, we say that $X$
is a $\mathbf{predecessor}$ of $Y$ or that $Y$ is a
$\mathbf{successor}$ of $X$, provided that there exists a path from
$X$ to $Y$. A tilting module $T$ is $\mathbf{convex}$ if, for a
given pair of indecomposable summands of $T$, $X$,$Y$  in add $T$,
any path from $X$ to $Y$ contains only indecomposable modules in add
$T$. Following \cite{APT}, we say that a set $\Sigma_{T}$ in mod
$A$ is a $\mathbf{complete}$ $\mathbf{slice}$ if
 $T = \bigoplus_{M \in \Sigma_{T}}M$ is a convex tilting module
with End$_A T$ hereditary. For the original definition of complete
slice, we refer the reader to \cite{Ri1}, \cite{Ri2}.

\vskip .1cm
For a given tilting module $T$ in mod $H$ there exists two
full disjoint subcategories of mod $H$, namely

$$\mathcal{F}(T)=\{X \in \mbox{mod}H \mbox{ such that } \mbox{Hom}_H(T,X)=0 \}$$
$$\mathcal{T}(T)=\{X \in \mbox{mod}H \mbox{ such that } \mbox{Ext}^{1}_H(T,X)=0 \}$$
\vskip .1cm

\noindent the free torsion class and  the torsion class, respectively.
\vskip .1cm
 Furthermore, if $T$ is a convex tilting module, then
mod$H=\mathcal{F}(T)\bigcup \mathcal{T}(T)$. We have that, in this
case,  $\mathcal{F}(T)$ is closed under predecessors and
$\mathcal{T}(T)$ is closed under successors.

\vspace{0.2in}

We will say that the category  $\mathcal{F}(T)$($\mathcal{T}(T)$) is finite if and only if there exists a finite number of iso-classes of indecomposable modules in $\mathcal{F}(T)$($\mathcal{T}(T)$), otherwise we say that the category $\mathcal{F}(T)$($\mathcal{T}(T)$) is infinite.

\vspace{0.2in}

Let us recall the following result from Happel and Ringel about tilting modules in tame hereditary algebras.

\vspace{0.2in}

\begin{lem}\cite[Lema 3.1]{HR}
Let $H$ a tame hereditary algebra. If $T$ is a tilting module over $H$ then $T$ has at least one non zero preprojective or preinjective direct summand.
\end{lem}

\vspace{0.2in}

The previous lemma implies that a tilting module over an hereditary tame algebra is not regular.
In fact, if we assume that all non regular summands of $T$ are preinjectives we obtain the following description of the categories $\mathcal{T}(T)$ and $\mathcal{F}(T)$.

\vspace{0.2in}

\begin{prop}\cite[Proposici\'on (3.2)*]{HR}\label{3.2}

Let $H$ be a tame hereditary algebra, and $T$ a tilting module over $H$. Then the following are equivalent:

\begin{itemize}
\item [(a)] $\mathcal{F}(T)$ is infinite.
\item [(b)] $\mathcal{F}(T)$ contains infinite preinjective indecomposable modules.
\item [(c)] There exists an indecomposable homogeneous simple module which belongs to $\mathcal{F}(T)$.
\item [(d)] All homogeneous modules belong to $\mathcal{F}(T)$.
\item [(e)] All preprojective modules belong to $\mathcal{F}(T)$.
\item [(f)] $T$ has no non zero preprojective summand.
\end{itemize}
\end{prop}

\vspace{0.2in}

Dualy, we have a similar description for a tilting module without preinjective direct summands.

\begin{prop}\cite[Proposici\'on 3.2]{HR}\label{3.2*}

Let $H$ be a tame hereditary algebra, and $T$ a tilting module over $H$. Then the following are equivalent:

\begin{itemize}
\item [(a)] $\mathcal{T}(T)$ is infinite.
\item [(b)] $\mathcal{T}(T)$ contains infinite preprojective indecomposable modules.
\item [(c)] There exists an indecomposable homogeneous simple module which belongs to $\mathcal{T}(T)$.
\item [(d)] All homogeneous modules belong to $\mathcal{T}(T)$.
\item [(e)] All preprojective modules belong to $\mathcal{T}(T)$.
\item [(f)] $T$ has no non zero preinjective summand.
\end{itemize}
\end{prop}

The following remark follows from (a) if and only if (f) and is very helpful.

\begin{obs}\label{remark de 3.2}
If $T$ has  preinjective direct summands then  $\mathcal{T}(T)$ is finite. Moreover, if the only non regular direct summands of $T$ are preinjective then $\mathcal{F}(T)$ is infinite.
\end{obs}

\vspace{0.2in}

Analogously, if a tilting module $T$ over a tame hereditary algebra has no non zero preinjective summand the roles of $\mathcal{F}(T)$ and $\mathcal{T}(T)$ are inverted, it is, $\mathcal{F}(T)$ is finite and $\mathcal{T}(T)$ is infinite.

\vspace{0.2in}

The following proposition allow us to determine when to expect that the category $\mathcal{F}(T)$ be finite.

\begin{prop}\cite[Proposici\'on 3.3]{AK}\label{salvaje}

Let $H$ be a representation infinite hereditary algebra and $(\mathcal{T},\mathcal{F})$ a torsion pair over mod $H$. Then $\mathcal{T}$ contains infinite classes of isomorphisms of indecomposable preprojective modules if and only if $(\mathcal{T},\mathcal{F})$ is induced by a tilting module without preinjective summands and $\mathcal{F}$ is finite.

Moreover, if $H$ is wild, then $(\mathcal{T},\mathcal{F})$ is induced by a preprojective tilting module.
\end{prop}

\vspace{0.2in}
Combining the last two results we obtain under what hypothesis a tilting module $T$ induces a torsion pair  $(\mathcal{T}(T),\mathcal{F}(T))$ whit one of the categories finite.

\vspace{0.2in}

For every indecomposable regular module we have the cone determined by $M$, defined as follows.

\begin{Def}\label{wing}\cite[Definici\'on XVI. 1.4]{SS2}
Given an indecomposable $H$-module  $M$ in a regular component $\mathfrak{C}$ of the Auslander-Reiten quiver $\Gamma_H$ of an hereditary algebra $H$, the cone determined by $M$ is a full subquiver of $\mathfrak{C}$ whith the following shape:

$$\xymatrix@!0 @R=10mm@C=10mm{
M_{1,1}\ar[rd] &&M_{2,2}\ar[rd]&&M_{3,3}\ar[rd]&&\cdots&& M_{m-2,m-2}\ar[rd]&&M_{m-1,m-1}\ar[rd]&&M_{m,m} \\
&M_{1,2}\ar[rd]\ar[ru]&&M_{2,3}\ar[rd]\ar[ru]&&\cdots&&\cdots\ar[ru]&&M_{m-2,m-1}\ar[rd]\ar[ru]&&M_{m-1,m}\ar[ru]& \\
&&M_{1,3}\ar[rd]&&\cdots&&\cdots&&\cdots\ar[ru]&&M_{m-2,m}\ar[ru]&& \\
&&&\cdots&&\cdots\ar[rd]&&\cdots&&\cdots\ar[ru]&&& \\
&&&&\cdots\ar[rd]&&M_{2,m-1}\ar[ru]\ar[rd]&&\cdots&&&& \\
&&&&&M_{1,m-1}\ar[rd]\ar[ru]&&M_{2,m}\ar[ru]&&&&& \\
&&&&&&M\ar[ru]&&&&&& }$$

\end{Def}

We denote by $\mathcal{C}(M)$ to the cone determined by $M$. We say that $M$ is the vertex of $\mathcal{C}(M)$ and that the cone is of level $m$.

\vspace{0.1in}

The next lemma is a bit technical but is really useful to deal with regular summands in a stable tube of a tilting $T$.

\begin{lem}\cite[XVII Lema 1.7 ]{SS2}\label{union de conos}

Let $H$ be an hereditary algebra and $\mathcal{R}$ a stable tube in $\Gamma_{H}$. For any pair of indecomposable modules $M,N$ in $\mathcal{R}$ such that:

$$\mbox{Ext}^{1}_H(M\oplus N, M\oplus N)=0$$

\noindent one of the following conditions is satisfied:

\begin{itemize}
\item [a)] $\mathcal{C}(M) \subset \mathcal{C}(N)$,
\item [b)] $\mathcal{C}(N) \subset \mathcal{C}(M)$ or
\item [c)] $\mathcal{C}(M) \bigcap \mathcal{C}(N)= \emptyset$, $\mathcal{C}(M) \bigcap \tau\mathcal{C}( N)= \emptyset$, and $\tau \mathcal{C}(M) \bigcap \mathcal{C}( N)= \emptyset$.
\end{itemize}

\end{lem}

\vspace{0.2in}

Let $T$ be a tilting module with non zero regular summands over a tame hereditary algebra $H$. Let $T_1,...T_n$ the regular summands of $T$ which belongs to a stable tube $\mathcal{R}$ of $\Gamma_H$. Since Ext$_{H}^{1}(T,T)=0$ by the previous lemma we can assume that $T_1,...,T_n$ belongs to a disjoint union of cones, $\bigcup_{i=1}^{r}C(T'_i)$ where each one of the modules $T'_i$ is a direct summand of $T$. Note that in general $r \leq n$. We will show later in section 4 that, under certain hypothesis, all the indecomposable regular modules in $\mathcal{T}(T)\bigcap\mathcal{R} $ belongs to the union of this cones.

\vspace{0.2in}

Every one of the cones recently described have in addition the following property:

\begin{prop}\cite[ XVII Proposici\'on  2.1]{SS}
Let $H$ be a hereditary algebra and $M$ an indecomposable module over a stable tube $\mathcal{R}$ such that $\mathcal{C}(M)$ is a cone of depth $m$. If $T$ is a tilting module such that $M$ is a direct summand of $T$ then $T$ have exactly $m$ indecomposable direct summands in $\mathcal{C}(M)$.
\end{prop}

\vspace{0.2in}

We will illustrate the results of this section with an example.

\vspace{0.2in}

\begin{ej}
Let $H=KQ$ the hereditary algebra given by the quiver $Q=\tilde{\mathbb{D}}_{12}$.

\vspace{0.2in}
$$\xymatrix @!0 @R=1.1cm  @C=1.1cm
{1\ar[rd]&&&&&&&&&&12\\
&3\ar[r]&4\ar[r]&5\ar[r]&6\ar[r]&7\ar[r]&8\ar[r]&9\ar[r]&10\ar[r]&11\ar[rd]\ar[ru]&\\
2\ar[ru]&&&&&&&&&&13}$$

Since $Q$ is an euclidean quiver , we know that the algebra $H$ is a infinite representation algebra of tame type.  Therefore we have that the Auslander-Reiten quiver of $H$, $\Gamma_H$ has a preprojective component, a preinjective component and infinite regular component which consists in a family of orthogonal stable tubes. In fact, in this example we have a tube $\mathcal{R}$ of rank 10. The mouth of this tube consists in nine simple modules that we will denote by  $E_i=S_{i+2}$ for $i=1,...,9$ and the module given by the representation:

\vspace{0.1in}
$$\xymatrix@!0 @R=1.1cm  @C=1.1cm
{&K\ar[rd]^{1}&&&&&&&&&&K\\
E_{10}:&&K\ar[r]^{1}&K\ar[r]^{1}&K\ar[r]^{1}&K\ar[r]^{1}&K\ar[r]^{1}&K\ar[r]^{1}&K\ar[r]^{1}&K\ar[r]^{1}&K\ar[ru]^{1}\ar[dr]_{1}&\\
&K\ar[ru]_{1}&&&&&&&&&&K}$$

\vspace{0.2in}
This modules also satisfy that $E_i=\tau E_{i+1}$ and $\tau E_1=E_{10}$.

\vspace{0.2in}

Let $T=\oplus^{12}_{i=1}T_i$ be the tilting module given by the following direct summands:

$$T_1=E_2^{1}, \, T_2= E_2^{3}, \, T_3=E_4^{1}, \, T_4=E_2^{4}, \, T_5=E_8^{1}, \, T_6=E_7^{2}, \, T_7=E^{1}_{10}$$

\vspace{0.1in}
\noindent where $E_{j}^{i}$ is the $i$-th module of the ray starting in $E_j$ and the other summands are the injective modules:

$$T_8=I_1, \, T_9= I_2, \, T_{10}=I_7, \, T_{11}=I_{10}, \, T_{12}=I_{12}, \, T_{13}=I_{13}.$$

\vspace{0.1in}
\noindent are the indecomposable direct preinjective summands of $T$.

\vspace{0.2in}

In this example, applying the proposition \ref{3.2*}, we get that the torsion $\mathcal{T}(T)$ is finite . We will only show how the torsion of $T$ looks like over the tube $\mathcal{R}$.

\vspace{0.1in}

 The part of the Auslander-Reiten quiver corresponding to $\mathcal{R}$ is given by the following translation quiver where the upper vertices represents the  modules $E_i$ in the mouth of $\mathcal{R}$ ordered from lowest to highest respect the index $i$:

\centerline{
\includegraphics[width=400pt]{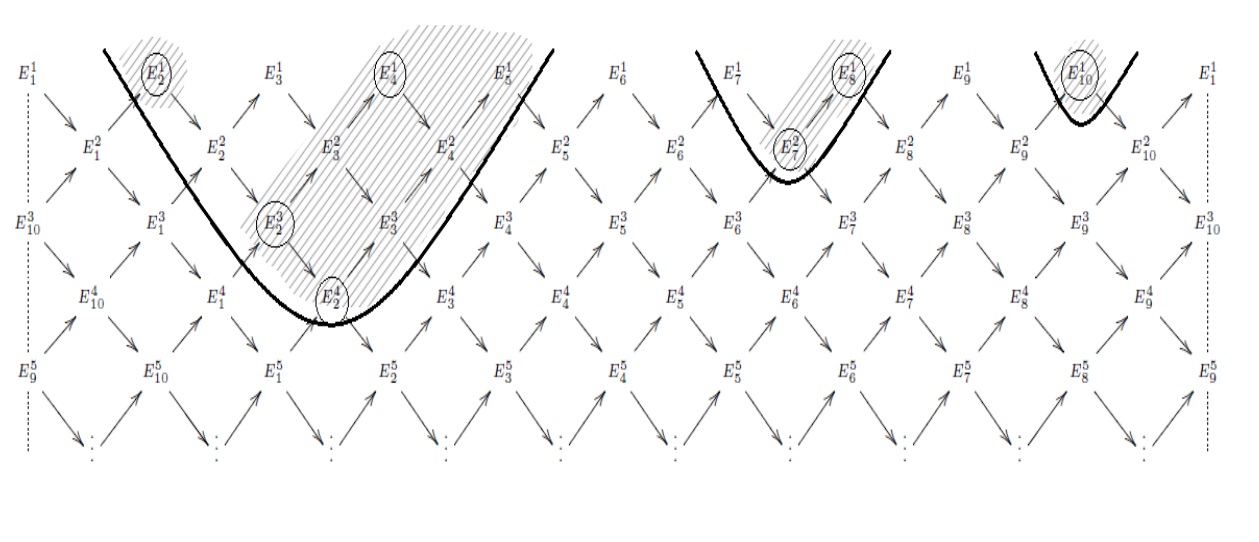}}

Here the  vertices within the  circles correspond to the summands of $T$. The shaded areas correspond to regular indecomposables modules in $\mathcal{T}(T)$. We can see that  the regular summands of $T$ are distributed over 3 cones . The modules $T_1,$ $T_2$ , $T_3$ and $ T_ 4$ belong to the cone $ \mathcal{C}(T_  4)$ of level 4. The modules $ T_ { 5 } $ and $ T_ { 6 } $ are in  the cone of level 2 determined by $ T_ { 6 } $ and finally $ T_ { 7 } $ belongs to the level 1 cone generated by itself . Furthermore, using  the mesh relations of the  Auslander - Reiten quiver on the component $\mathcal{R}$ and the fact that the rest of the summands of $T$ are preinjective it can be shown that the rays generated by the regular simple modules $ S_3 = E_1 $ , $ S_4 = E_6 $ and $ S_7 = E_9 $ are entirely contained in $\mathcal{F}(T)$. Somehow these rays separate the cones.

\vspace{0.2in}

The rest of the  indecomposable modules  in $\mathcal{T}(T)$ corresponds to preinjective modules. Moreover, due to the choice of $T$ all the preprojective modules are in
$\mathcal{F} (T)$, as well as all regular modules that do not belong to the tube $\mathcal{R}$.

\end{ej}

\subsection{Cluster categories and cluster tilted Algebras}

For the convenience of the reader, we start this section recalling
some definitions and results of cluster categories from
\cite{BMRRT}. Let $\mathcal{C}$ be the  $\mathbf{cluster}$
$\mathbf{category}$ associated to $H$ given by
$\mathcal{D}^{b}(H)/F$, where $F$ is the composition functor
$\tau_{\mathcal{D}}^{-1}[1]$. We represent by $\tilde{X}$ the class
of an object $X$ of $\mathcal{D}^{b}(H)$ in the cluster category. We
recall that Hom$_{\mathcal{C}}(\tilde{X},\tilde{Y})= \bigoplus_{i
\in \mathbb{Z}}\mbox{Hom}_{\mathcal{D}^b(H)}(X,F^iY).$ We also recall
that $S = \mbox{ind}H \bigcup H[1]$ is a fundamental domain of
$\mathcal{C}$. If $X$ and $Y$ are objects in the fundamental domain,
then we have that $\mbox{Hom}_{\mathcal{D}^b(H)}(X,F^iY)=0$ for all
$i \neq 0,1$. Moreover, any indecomposable object in $\mathcal{C}$ is of the form
$\tilde{X}$ with $X \in S$. \vskip .1cm

We say that $\tilde{T}$ in $\mathcal{C}$ is a $\mathbf{tilting}$
$\mathbf{object}$ if Ext$^{1} _{\mathcal{C}}(\tilde{T}, \tilde{T}) =
0$ and $\tilde{T}$ has a maximal number of non-isomorphic direct
summands. A tilting object in $\mathcal{C}$ has finite summands.
\vskip .1cm

There exists the following nice correspondence between tilting
modules and basic tilting objects. \vskip .1cm
 \noindent$
\mathbf{Theorem}$ \cite[Theorem 3.3.]{BMRRT} {\it\begin{itemize}
\item[(a)] Let T be a basic tilting object in $\mathcal{C} = \mathcal{D}^b(H)/F$, where $H$
is a hereditary algebra with $n$ simple modules.
\begin{itemize}
\item[(i)] $T$ is induced by a basic tilting module over a hereditary algebra $H'$,
derived equivalent to $H$.
\item[(ii)] $T$ has $n$ indecomposable direct summands.
\end{itemize}
\item[(b)] Any basic tilting module over a hereditary algebra $H$ induces a basic tilting
object for $\mathcal{C} = \mathcal{D}^b(H)/F$.
\end{itemize}}
\vskip .1cm
 In \cite{BMR2}, Buan, Marsh and Reiten introduced the
cluster-tilted algebras as follows. Let $\tilde{T}$ be a tilting object over the
cluster category $\mathcal{C}$. We recall that $B$ is a
$\mathbf{cluster}$ $\mathbf{tilted}$ $\mathbf{algebra}$ if $B=$
End$_{\mathcal{C}}(\tilde{T})$. It is also known that, if $\tilde{T}$ is a tilting object in $\mathcal{C}$, the
functor $\Hom \,\,)$ induces an equivalence of categories between
$\mathcal{C}/\add$ and mod $B$. We will call this equivalence the [BMR]-equivalence.

Thus using the above equivalence, we can compute Hom$_B(X',Y')$
in terms of the cluster category $\mathcal{C}$ as follows:
$$Hom_{B}(\Hom \tilde{X}),\Hom \tilde{Y}))\simeq \mbox{Hom}_{\mathcal{C}}(\tilde{X},\tilde{Y})/\add,$$ where $X'= \Hom \tilde{X})$ and $Y'= \Hom \tilde{Y})$ are $B$-modules.

\begin{obs} \label{2.4}
Suppose $X' = \Hom \tilde{X})$ and $Y'= \Hom \tilde{Y})$. If for every $f:\tilde{X}
\rightarrow \tilde{Y}$, $f$ factors through $\add$ in $\mathcal{C}$ then we have that $Hom_B(X',Y')=0$.
\end{obs}

\section{ Resolutions for tame cluster tilted algebras}


\subsection{Note on finite torsion class hereditary algebras}

In this section we will study the \rd of cluster tilted algebras given by a tilting module with finite torsion class. It follows from proposition \ref{salvaje} that if $H$ is an hereditary algebra of infinite representation type and $T$ is a tilting module over $H$ there is only two cases  where the torsion class is finite. Either the algebra $H$ is tame and $T$ has no nonzero preprojectives summands or $H$ is wild or tame and $T$ is a preinjective tilting $H$-module. The cluster tilted algebras arising from the second case hypothesis are exactly the  cluster concealed algebras. We known by \cite{GT} that such class of algebras has representation dimension three. Thus in this paper we will consider only the first case. We will show that the weak representation dimension of tame cluster tilted algebras is at most three.

\subsection{Tame cluster tilted algebras}

Let $B$ a tame cluster tilted algebra. Then we  have that $B=$End$_{\mathcal{C}}(\tilde{T})$ where $\tilde{T}$ is a cluster tilted object induced by a tilting module $T$ over an hereditary algebra $H$ and $\tilde{C}$ is the cluster category asociated to $H$. We know that $B$ is tame if and only if $H$ is tame, therefore by 
\cite{HR} $T$ has at least a non zero non regular summand. Moreover, due to the construction of cluster categories we can assume that all the non zero non regular indecomposables summands of $T$ are either preprojective or preinjective. Thus by \cite[Proposici\'on (3.2),(3.2)*]{HR} we can assume that either $\mathcal{F}(T)$ is finite or $\mathcal{T}(T)$ is finite respectively.

\vspace{0.2in}

In fact, if $\mathcal{T}$ is a cluster tilting object induced by an $H$-tilting module $T$ with non zero preinjective and non zero regular summands such that $\mathcal{T}(T)$ is finite, then there exists a hereditary algebra $H'$ derived equivalent to $H$ and a $H'$-module $T'$ with at least a non zero regular summand such that $\mathcal{F}(T')$ is finite in mod$H'$. Even more, the cluster category induced by $H$ and $H'$ is the same and clearly  $T'$ induces the same cluster tilted object $\tilde{T}$ over $\mathcal{C}$. We will ilustrate this situation in an example.

\vspace{0.2in}


\begin{ej}\label{ej2}
Let $H$ the algebra given by the dynkin diagram $\tilde{\mathbb{D}}_5$. We consider the tilting module $T = \tau^{-2}P_5 \oplus \tau^{-2}P_6 \oplus \tau^{-2}P_4 \oplus \tau^{-2}P_1 \oplus \tau^{-2}P_2 \oplus S_3$. This module has 4 indecomposable preprojective summands and the indecomposable regular summand $S_3$.

\vspace{0.2in}

Let $B =$End$_\mathcal{C}(\tilde{T})$. $B$ is giben by the following quiver:

\vspace{0.1in}
\begin{displaymath}
\def\objectstyle{\scriptstyle}
\def\labelstyle{\scriptstyle}
\xymatrix @!0 @C=16mm
{& 3 \ar[dl]_{\alpha}\ar[dr]^{\beta}& \\ 1 \ar[rd]_{\gamma}& & 2\ar[ld]^{\delta}
 \\ & 4 \ar[rd]_{\lambda}\ar[ld]^{\mu}\ar[uu]^{\epsilon} &
\\5& &6
}
 \end{displaymath}

\vspace{0.1in}
\noindent with the relations $\epsilon\alpha=\epsilon\beta=\gamma\epsilon=\delta\epsilon=0$ and $\alpha\gamma=\beta\delta$.
\end{ej}

\vspace{0.2in}
In this example we have that the category  $\mathcal{F}(T)$  is finite since $T$ has no non zero preinjective summands.
\vspace{0.2in}

 Now consider the following tilting $H$-module  $T' = \tau^{2}I_5 \oplus \tau^{2}I_6 \oplus \tau^{2}I_4 \oplus \tau^{2}I_1 \oplus \tau^{2}I_2 \oplus S_3$. Unlike the previous module, here all the indecomposable summands are preinjective to exception  of the regular summand $S_3$ and  the torsi \ ' on associated with this module is finite.

\vspace{0.2in}

Finally, we can see that this module induces the same  cluster tilted algebra $B=$ End$_{C}(\tilde{T})\simeq$ End$_{\mathcal{C}}(\tilde{T}')$ and therefore the algebra $B$ can be studied beginning with any of the $H$-modules $T$ or $T'$ . As shown in this example here the algebras $H$ and $H'$ coincide. In general this is not necessarily true. One must consider all possible changes in orientation on the ordinary quiver  of $H$ in order to obtain an adequate algebra $H'$. Also note that since we only consider changes of orientation the algebra $H'$ is always derived equivalent to $H$.

\vspace{0.2in}


\vspace{0.2in}

From now on, we will assume that any tame cluster tilted algebra End$_{\mathcal{C}}(\tilde{T})$ is induced by a tilting module $T$ such that every indecomposable summand of $T$ is either preinjective or regular.  Thus by  lemma \ref{3.2*} we know that the category $\mathcal{T}(T)$. We write $T = T_{I} \oplus T_{R}$ where $T_{I}$ is a preinjective module and $T_{R}$ is a regular module. We recall that under these hypotheses $T_{I}$ must be a nonzero summand of $T$. We will only consider the case when $T_{R}\neq 0$ as we discussed before.

\vspace{0.2in}

Since we work under the assumption that $\mathcal{T}(T)$ is finite we can define the following $H$-module.
Let $\mathcal{T}$ the direct sum of all indecomposable non isomorphic modules in $\mathcal{T}(T)$. Thus applying the hom$_{\mathcal{C}}(\tilde{T}, \;)$ to $\mathcal{T}$ we obtain the $B$-module $\mathcal{T}'=\Hom \tilde{\mathcal{T}})$. We will show that we only need the modules $\mathcal{T}'$ and $H'$ to get an approximation of all, except for finitely many, indecomposable modules in mod $B$.

\vspace{0.2in}

Since our objective is to compute the representation dimension of $B$, we need to construct a generator module $M'$ for mod $B$.

\vspace{0.2in}

We start constructing a transjective module $M'_1$. First, let $\mathcal{Y}$ be the the direct sum of all indecomposable preinjective modules, up to isomorphisms, in $\mathcal{T}$. Recall that if $\Sigma$ is a tilting module given by a complete slice then any indecomposable module belongs to $\mathcal{T}(\Sigma)$ if and only if is a  successor of $\Sigma$. Since $\mathcal{T}(T)$ is finite, it is possible to choose a preinjective slice module $\Sigma$ such that $\mathcal{Y}$ belongs a $\mathcal{T}(\Sigma)$. Since $\Sigma$ is preinjective then $\mathcal{T}(\Sigma)$ is also finite. Hence we can define the following module.

\vspace{0.2in}

Let $N_1=\oplus_{X \in \mathcal{T}(\Sigma)}X$ with $X$ indecomposable. This module induces the class $\tilde{N_1}$ in the cluster category $\mathcal{C}$. We take the image of this object under the functor hom$_{\mathcal{C}}(\tilde{T}, \;)$, which we will denote by $N'_1$. Finally we set $M'_1 = N'_1 \oplus \mathcal{Q}' \oplus H'$. Here $\mathcal{Q}'= \bigoplus Q'$ tal que $Q' \simeq \Hom \tilde{Q})$ with $Q = P[1]$ where $P$ is an indecomposable projective module in mod $H$ and $H'= \Hom \tilde{H})$ where  $H$  represents the algebra as a projective $H$-module.

\vspace{0.2in}

 We have the following description for $M'_1$. Let $H_1$ be the hereditary algebra $\tau_{\mathcal{D}}^{-2}H$. Then, the $B$-module $M'_1$ is obtained by applying the functor hom$_{\mathcal{C}}(\tilde{T}, \;)$ to the object $\tilde{M_1}$ in $\mathcal{C_{H_{1}}}$. Here $M_1$ is the $H_1$-module given by $\oplus_{X \in \mathcal{T}(\Sigma)}X$ where $\Sigma$ is the same complete slice choosen for mod $H$ but instead considered as $H_1$-module and $\mathcal{T}_{H_1}(\Sigma) \subset $ mod $H_1$. Note that the object $M_1$ in $D^{b}(H)$ is obtained as the direct sum of all indecomposable objects which are both successors of $\Sigma$ and predecessors of the slice $\tau_{\mathcal{D}}^{-1}H[1]$. This mean that the $B$-module $M'_1$ is between  two local slices (see \cite{ABS}) in mod $B$. This way to visualize the module $M'_1$ is not only more simple but will also help to simplify some technical details in some proofs.

\vspace{0.2in}

Before we complete the definition of the module $M'$,  we need to establish the following definition:

\vspace{0.2in}

\begin{Def}
Let $E$ be a regular module in a stable tube $\mathcal{T}_{\lambda}$.
We will call \textbf{edge of the cone} $\mathcal{C}(E)$ and denoted by $\mathcal{C}(E \to)$ to the intersection of the coray containing  $E$ in $\mathcal{T}_{\lambda}$ and the cone $\mathcal{C}(E)$.
\end{Def}

Note that  the edge of the cone cogenerates all modules in $\mathcal{C}(E)$ and that $E$ generates all modules in the edge of the cone.

\vspace{0.2in}

Let $T_{1}$,...,$T_{n}$  be all the indecomposables regular summands of $T_R$. For each $T_i$ we set $\mathcal{C}(T_{i}),$ $  i=1,...,n$ the cone of $T_i$. We have the following observations:

\begin{itemize}

\item [i)] Every summand $T_{1}$,...,$T_{n}$ belong to at least one of the cones $\mathcal{C}(T_i)$ but possible more than one.
\item [ii)] Every one of the vertex of this \textbf{cones} is obviously a indecomposable summand of $T_{R}$.
\item[iii)] The edge of the cone $\mathcal{C}(T_i \rightarrow)$ is enterally contained in $\mathcal{T}(T)$ for every $i=1,...,n$.
\end{itemize}

\vspace{0.1in}

We define $W_i$ as the regular $H$-module given by $\bigoplus_{X \in \mbox{ind}(\mathcal{C}(T_i))} X $.

\vspace{0.1in}

Since we are assuming that $T_R$ has $n$ regular indecomposable summands, then we have $r$ modules $W_1$,...,$W_r$  with the above properties. However , we can choose by lemma \ref{union de conos} regular summands of $ T $, ( rearranging  the index of the summands) $T_1, ..., T_m$ such that we can replace  i) by:

\begin{itemize}
\item [i')] Every one of the summands $T_{1}$,...,$T_{n}$ belongs to one and only one of the cones $\mathcal{C}(T_j)$ for $j=1,...,m$.
\end{itemize}

\vspace{0.1in}
We call this cones \textbf{maximal} respect to property  i). More over, it is possible to determine the exact amount of  indecomposable summands  of $T$ in each of this cones by only  knowing  the level of the  vertex $T_j$ in the stable tube.

\vspace{0.2in}

Let $M'_2$ be $(W_1)' \oplus ... \oplus (W_m)'$ where each $(W_j)'$ denotes the image of the object $\tilde{W}_j$ by the  equivalence hom$_{\mathcal{C}}(\tilde{T}, \;)$. Note that this module is not transjective neither regular in mod $B$.

\vspace{0.2in}

Finally, the $B$-module $M'$ is given by $M' = M'_1 \oplus M'_2$. We will show that $M'$ satisfies (a) of definition \ref{repdim} for $i=3$.

\subsection{The approximation resolution}

The aim of this subsection is to study when is possible to obtain an exact sequence:

\vspace{0.1in}
$$0\rightarrow M_{1} \rightarrow M_{0} \stackrel{f'}\rightarrow X' \rightarrow 0$$

\vspace{0.1in}
\noindent for each $X'$ in ind $B$ where $M_{i} \in $ add$(M')$ for $i=0,1$.

\vspace{0.2in}

If $X' \in$ add$(M')$, the existence of such sequence is trivial. Since $\mathcal{Q}'$ is in add$(M')$, we will assume that $X'=\Hom \tilde{X})$ with $X \in $ ind $H$, it is, $X'$ will always be obtained from a $H$-module $X$.

\vspace{0.2in}

Let $X$ be in ind $H$.  Associated to $X$ we have the modules $tX$, $X/tX$  and the exact sequence given by

$$0  \rightarrow   tX  \stackrel{i}\rightarrow   X  \stackrel{p}\rightarrow   X/tX  \rightarrow  0.$$

We will construct a new short exact sequence  from a minimal projective resolution of  $X/tX$.

\vspace{0.2in}

Let

$$0\rightarrow K \rightarrow P_{X/tX} \stackrel{\pi}\rightarrow X/tX \rightarrow 0$$

be a minimal projective resolution of $X/tX$ where $\pi:P_{X/tX} \rightarrow X/tX$ is a minimal projective cover of $X/tX$.

\vspace{0.2in}

Since $p:X \to X/tX$ is an epimorphism, there exists a morphism $g :P_{X/tX} \to X$ such that $pg=\pi$. Hence, we set
$f:=\left(\begin{array}{cc}i,\,  g\end{array}\right)$ and  we have the following commutative diagram:

\vspace{0.1in}
$$
\xymatrix{
&0&&0&\\
0\ar[r]& tX\ar[u] \ar[r]^{i} & X \ar[r]^{p} & X/tX\ar[u] \ar[r] & 0 \\
0 \ar[r] & tX\ar[u]^{Id} \ar[r] & tX \oplus P_{X/tX}\ar[u]^{f} \ar[r] & P_{X/tX}\ar[u]^{\pi} \ar[r] & 0 \\
&   0\ar[u] &   &   K\ar[u]   &  \\
&  &    &   0 \ar[u]   &
}
$$

\vspace{0.1in}
\noindent More over, is easy to check that the columns and rows are exact sequences.

\vspace{0.2in}

Applying the snake lemma, we obtain a short exact sequence:

\begin{equation}\label{ec 1}
 0 \to K \to tX \oplus P_{X/tX} \stackrel{f}\to X \to 0 .
\end{equation}
\vspace{0.2in}

Our objective is to study when this short exact sequences remains exact in mod $B$. It is, we will show under which hypothesis this short exact sequence induces a short exact sequence

\vspace{0.2in}
$$0\rightarrow K' \rightarrow  (tX)' \oplus (P_{X/tX})' \stackrel{f'}\rightarrow X' \rightarrow 0$$

\vspace{0.1in}
\noindent in mod $B$.

\vspace{0.1in}
\begin{prop}\label{resolucion}
Let $X'=$Hom$_{\mathcal{C}}(\tilde{T},\tilde{X})$ with $X \in$  mod $H$.
Then the morphism $f:tX \oplus P_{X/tX} \stackrel{f}\to X$ induces an epimorphism $(tX)' \oplus P'_{X/tX} \to X'$. More over, if $X \notin \mathcal{T}(T)$ then the  sequence
$$0 \to K' \to (tX)' \oplus P'_{X/tX} \to X' \to 0$$
\noindent is exact in mod$B$.
\end{prop}

\vspace{0.2in}

Before we proceed with the proof of proposition \ref{resolucion},  we will discuss some details that show up in the proof.

\vspace{0.2in}

First, we need to assure when a morphism in $\mathcal{C}$ induces a non zero morphism in mod $B$. Recall that

\begin{obs} \label{null morphism}
If $X' = \Hom \tilde{X})$ and $Y'= \Hom \tilde{Y})$ are indecomposable in mod $B$, then a morphism $f'$ is non zero if the morphism  $\tilde{f}:\tilde{X} \rightarrow \tilde{Y}$ does not factorizes  through $\add$ in $\mathcal{C}$.
\end{obs}

\vspace{0.2in}











\vspace{0.2in}

The next theorem gives us necessary and sufficient  condition to determine when the morphism $h: \tilde{X} \to \tilde{E[1]}$ factors through  add$(\tau \tilde{T})$ in $\mathcal{C}$.

\vspace{0.2in}

\begin{teo}\label{fact tau T}
Let $\Sigma$ be a preprojective complete slice in mod $H$ such that $\tau T \in \mathcal{T}(\tau^{-1}\Sigma)$.
Let $X \in \mathcal{T}(\tau^{-1}\Sigma)$ and let $E$ be a $H$-module in add $(\Sigma)$.

Let
$$E \to M \to X \stackrel{h}\to E[1]$$ be a triangle in $D^{b}(H)$.

Therefore the morphism $h:\tilde{X} \to \tilde{E[1]}$ factors through add$(\tau \tilde{T})$ in $\mathcal{C}$ if and  only if exists a morphism $k:X \to T'$ with $T' \in $ add$(\tau \tilde{T})$ such that $hi=0$ where $i$ is the canonic inclusion of Ker $k$ in $X$.
\end{teo}

\vspace{0.1in}
\noindent \textbf{proof:}

Let $X$ be an $H$-module that does not belong to add$(\Sigma)$ and assume that we have a triangle

\vspace{0.1in}
$$ E \to M \to X \stackrel{h}\to E[1].$$

\vspace{0.2in}
Proceeding as in \cite[proposition 3.2]{GT}, let $H_1 =\tau^{-1}_{\mathcal{D}}\Sigma$, is easy to see that $X, \tau T, T$ can be identified with $H_1$-modules and that the slice $\Sigma[1]$ can be identified  with the complete slice formed by the indecomposable injective $H_1$-modules.

\vspace{0.2in}

Suppose that we have a non zero morphism $k:X \to T'$ con $T' \in$ add$(\tau T)$. Let $K$ be the kernel of $k$, Ker $k$ and let $i$ be the inclusion of kernel of $k$ in $X$, then we have that

\vspace{0.1in}
 $$0 \to K \stackrel{i} \to X \stackrel{k} \to T' $$

\vspace{0.1in}
\noindent is exact in mod $H_1$. We have that $E[1]$ is an injective module  considered as $H_1$-module, then the sequence:

\vspace{0.2in}
 $$ \mbox{Hom}_{H_1}(T',E[1]) \stackrel{k^{*}}  \to \mbox{Hom}_{H_1}(X,E[1])  \stackrel{i^{*}} \to \mbox{Hom}_{H_1}(Ker k,E[1]) \to 0$$

\vspace{0.1in}
\noindent is exact in mod $H_1$, where $k^{*}= $ Hom$_{H_1}(k, E[1])$ e $i^{*}=$ Hom$_{H_1}(i,E[1])$. Since $Im k^{*} = $ Ker $i^{*}$,  $hi=0$ if and only if  exists a non zero morphism, $l:T' \to E[1]$ such that $h=kl$. Therefore, we have the following diagram:

\vspace{0.1in}
\begin{displaymath}
\def\objectstyle{\scriptstyle}
\def\labelstyle{\scriptstyle}
\xymatrix @!0 @C=16mm
{  &  & \mbox{Ker}k\ar[d]^{i} & \\
K\ar[r] & E\ar[r] & X\ar[d]^{k}\ar[r]^{h} & E[1]\\
&& T'\ar[ru]^{l}&}
 \end{displaymath}

\vspace{0.1in}
\noindent Hence, $h$ factors through add$(\tau T)$ if and  only if $hi=0$.$\Box$


\vspace{0.2in}

This gives us a necessary and sufficient condition  under which the exact sequence \ref{ec 1} induces a short exact sequence such as the described in proposition \ref{fact tau T}.

\vspace{0.2in}

\begin{obs}\label{NOTA}
Recall from \cite{GT}, that an epimorphism $f:M \to N$ in mod $H$ induces an epimorphism in mod $B$ if and only if all morphism $\tilde{h}:\tilde{N} \to \tau $ker$\tilde{f}$ factors through $\add$ in the cluster category $\mathcal{C}$.
\end{obs}

\vspace{0.2in}

Therefore, as an immediate consequence of the previous theorem we can conclude that any triangle of the form $K \to E \to X \to K[1]$ satisfying  the condition described in \ref{NOTA}, induces an epimorphism $E' \to X'$ in mod $B$. We may establish a dual version of this theorem which gives  conditions necessary and sufficient for a triangle induces a monomorphism in mod $B$.

\vspace{0.2in}

Note that the condition required in the Theorem\ref{fact tau T}, $hi=0$, it is satisfied trivially if $X'$ is of the form $\Hom X)$, with $X \in \mathcal{F}(T)$. In fact, since $\mathcal{F}(T)=$Cogen$(\tau T)$, exists a monomorphism $k:X \to T'$ with $T' \in $ add$(\tau T)$. Hence, Ker$k=0$ and therefore $i=0$. Which implies that $hi=0$ for every $h$. This shows that in this particular case, the sequence described in the Proposition\ref{resolucion} is similar to the described in \cite[proposition3.2]{GT}.

\vspace{0.2in}

We proceed with proof of proposition \ref{resolucion}.

\vspace{0.2in}

\noindent \textbf{Proof of proposition \ref{resolucion}:}

\vspace{0.1in}
Let $X'$ be  such that $X$ is a non projective indecomposable $H$-module which does not belong to $\mathcal{T}(T)$. Observe that if $X$ is projective we have that $X'\in $ add$(M')$ and there is nothing to prove.

\vspace{0.2in}

Recall that, in the beginning of this section, we had constructed an exact sequence in mod $H$ of the form:

$$ 0 \to K \to tX \oplus P_{X/tX} \stackrel{f}\to X \to 0 $$.

\vspace{0.2in}

This exact sequences induces a triangle in $\mathcal{D}^{b}(H):$

\vspace{0.1in}
$$K \rightarrow tX \oplus P_{X/tX} \rightarrow X \stackrel{h}\rightarrow K[1].$$

\vspace{0.2in}

We want to show that this triangle induces an exact sequence in mod $B$. We will prove that $f':(tX)' \oplus P'_{X/tX} \to X'$ is an epimorphism. As we state before it suffices to show that $h:\tilde{X} \to  \tau\tilde{K}$ factors through add$(\tau\tilde T)$ in $\mathcal{C}$.

\vspace{0.2in}

Since Cogen$(\tau T)=\mathcal{F}(T)$, and $X/tX \in \mathcal{F}(T)$  we know that there  exists a monomorphism $j:X/tX \to T'$ with $T' \in $add$(\tau T)$. Since we assume that $X \notin \mathcal{T}(T)$, $X/tX \neq 0$ and therefore $p:X \to X/tX$ is a non zero morphism. Hence, the composition $jp$ gives rise to a non zero morphism $k:X \to T'$. More over, note that Ker $k = tX$. let $i:tX \to X$ be the canonical inclusion. Recall that $K$ is the kernel of the projective cover of $X/tX$ and hence is a projective module. Let $\Sigma$ be the projective complete slice, by Theorem\ref{fact tau T}, it suffices to prove that the composition $hi:tX \to K[1]$ is zero.

\vspace{0.2in}

From the exact sequence  (\ref{ec 1}) , we have that $h\left(\begin{array}{cc}i, \,g\end{array}\right) = \left(\begin{array}{cc}hi, \,hg\end{array}\right) = 0$. Hence $hi=hg=0$.

\vspace{0.2in}

Therefore $f':(tX)' \oplus P'_{X/tX} \to X'$ is an epimorphism in mod $B$.

\vspace{0.2in}

We will show that this epimorphism can be completed to an exact sequence:

\vspace{0.1in}
$$0\to (K)' \to (tX)' \oplus P'_{X/tX} \to X' \to 0 $$

\vspace{0.1in}
\noindent with $K' \in $add$(M')$.

\vspace{0.2in}

First, we start with the minimal projective resolution of $X$ in mod $H$:

\vspace{0.2in}
$$0 \to K \stackrel{j}\to tX \oplus P_{X/tX} \to X \to 0$$

\vspace{0.2in}

This short exact sequence induces the following triangle:

\vspace{0.2in}

$$ \tilde{K} \stackrel{\tilde{j}}\to \tilde{tX} \oplus \tilde{P}_{X/tX} \to \tilde{X} $$

\vspace{0.2in}
\noindent in the cluster category $\mathcal{C}$.

\vspace{0.2in}
Applying the functor Hom$_{C}(\tilde{T}, $\;\;\;$ )$ to this triangle we obtain the following sequence in $\mathcal{C}$/add$(\tau \tilde{T})$:

\vspace{0.2in}

$$\mbox{Hom}_{C}(\tilde{T},\tilde{K} )\stackrel{\tilde{j}^{*}}\to \mbox{Hom}_{C}(\tilde{T},\tilde{tX} \oplus \tilde{P}_{X/tX}) \to ...$$

\vspace{0.2in}

We will show that the morphism $\tilde{j}^{*}$ given by composing with $\tilde{j}$, it is, $\tilde{j}^{*}=\mbox{Hom}_{C}(\tilde{T},\tilde{j} )$ is a monomorphism.

\vspace{0.2in}

Let $\tilde{h}\in \mbox{Hom}_{C}(\tilde{T},\tilde{K} )$ such that $\tilde{h} \in$ Ker $\tilde{j}^{*}$, it is, $\tilde{j}\tilde{h}=0$. Since $\tilde{h}\in \mbox{Hom}_{C}(\tilde{T},\tilde{K} )$ we know that $\tilde{h}=(h_1, h_2)$ with $h_1:T \to K$ y $h_2: T \to \tau^{-1}K[1]$. Since we choose $T$ without preprojective summands, then Hom$_{H}(T,K)=0$ and $h_1=0$. Hence $\tilde{h}=h_2: T \to \tau^{-1}K[1]$.

\vspace{0.2in}

Since $Im j \subset P_{X/tX}$, the composition $\tilde{j}\tilde{h}\in \mbox{Hom}(\tilde{T},\tilde{P}_{X/tX})$. More over,
the morphism $\tilde{j}\tilde{h}$ is induced by the morphism $Fjh_2: T \to \tau^{-}1[1] \to \tau^{-1}P_{X/tX}[1]$ in $\mathcal{D}^{b}(H)$, where $F= $Hom$_{C}(\tilde{T}, $\;\;\;$ )$

\vspace{0.2in}

By hypothesis $\tilde{j}\tilde{h}=0$ in $\mathcal{C}$/add$(\tau T)$, therefore  $Fjh_2: T \to \tau^{-}K[1] \to \tau^{-1}P_{X/tX}[1]$ must factor through $\tau T$ en $\mathcal{D}^{b}(H)$. But, since $T$ is a tilting module, Hom$_{H}(T, \tau T)=0$ and thus this factorization implies that  $Fjh_2=0.$

\vspace{0.2in}

We will prove that $Fj$ is a monomorphism in mod $\tau_{D}^{-2}H$.

\vspace{0.2in}

First, we will show that $\tau^{-1}j$ is a monomorphism. Let $E_{X}= tX \oplus P_{X/tX}$
Applying $D$ the dual functor to the short exact sequence (\ref{ec 1}) we obtain:

\vspace{0.2in}

$$0 \to DX \stackrel{j}\to DE_{X} \to DK \to 0$$

\vspace{0.2in}
Applying Hom$( $\;\;\;$, H)$ to the last sequence we have:

\vspace{0.2in}

$$\mbox{Hom}(DK, H) \to \mbox{Hom}(DE_X, H) \to \mbox{Hom}(DX,H) \to $$
$$\to \mbox{Ext}^{1}(DK,H) \to \mbox{Ext}^{1}(DE_X, H)$$

\vspace{0.2in}

Since $X$ has no injective summands (Recall that $X \notin \mathcal{T}(T)$) therefore $DX$ has no projective summands and hence $\mbox{Hom}(DX,H)=0$. This implies that $\mbox{Ext}^{1}(DK,H) \to \mbox{Ext}^{1}(DE_X, H)$ is a monomorphism and since $DK$ and $DE_X$ has no projective summands we conclude that $\mbox{Ext}^{1}(DK,H)= tr DK= \tau^{-1}K$ and $\mbox{Ext}^{1}(DE_X, H)= tr DE_X = \tau^{-1}E_X$. Hence $\tau^{-1}j: \tau^{-1}K \to \tau^{-1}E_X$ is a monomorphism and $Fj: \tau^{-}1K[1] \to \tau^{-1}P_{X/tX}[1] $ is a monomorphism in mod $\tau^{-2}_{\mathcal{D}}H$.

\vspace{0.2in}
Finally, we get that $h_2=0$ since $Fj$ is a monomorphism and $Fjh_2=0$. Therefore $\tilde{h}=0$ and $\tilde{j}^{*}$ is a monomorphism.

Thus we have proved that the sequence:

\vspace{0.2in}
$$0\to (K)' \to (tX)' \oplus P'_{X/tX} \to^{f} X' \to 0 $$

\noindent is exact in mod $B$ $\Box.$
\vspace{0.2in}

 Last proposition shows that the chosen module $M'$ satisfies the first part of property (a) on the definition \ref{repdim} for every $X' \in $ mod $B$  of the form $X' = \mbox{Hom}_{\mathcal{C}}(\tilde{T},\tilde{X})$ such that $X \notin \mathcal{T}(T)$. But, since we are only considering the case in that $\mathcal{T}(T)$ is finite, then we only have a finite number of indecomposable modules (up to isomorphism) which are not of this form. Let $\mathcal{G'}= \oplus_{\mbox{ind }B} X'$ such that $X' = \mbox{Hom}_{\mathcal{C}}(\tilde{T},\tilde{X})$ with $X \in \mathcal{T}(T)$.

\vspace{0.2in}

We will show that $\mathcal{G'}$ is in fact a direct summand of $M'$ and thus we have a short exact sequence:

\vspace{0.2in}
$$0\rightarrow M_{1} \rightarrow M_{0} \stackrel{f'}\rightarrow X' \rightarrow 0$$

\vspace{0.1in}
\noindent  with $M_{i} \in $ add$(M')$, $i=0,1$ for every $X' \notin $add$(M')$ in mod $B$.

\vspace{0.2in}

\begin{teo}\label{prop torsion}
Let $H$ be a tame hereditary algebra  and  $T$ a tilting $H$-module without preprojective direct summands.
Let $X$ be a non zero indecomposable regular module in $\mathcal{T}(T)$. Then there exists a indecomposable regular summand of  $T$, namely $E$, such that $X \in \mathcal{C}(E)$.
\end{teo}

\textbf{proof:}

Suppose $X$ is a non zero indecomposable regular module in $\mathcal{T}(T)$. Then since $H$ is tame there exists a stable tube $\mathcal{T}_{\lambda}$ such that $X \in \mathcal{T}_{\lambda}$. More over, since $T$ has no preprojective summands and any two regular components of the Auslander-Reiten quiver of $H$ are orthogonal, the fact that $X \in \mathcal{T}(T)$ implies that there is at least an indecomposable regular summand of $T$ in $\mathcal{T}_{\lambda}$. Otherwise, Hom$_{H}(T,\mathcal{T}_{\lambda})=0$, which is a contradiction to the fact that $X \in \mathcal{T}(T)\bigcap\mathcal{T}_{\lambda}$.

\vspace{0.2in}

Let $T_1,...,T_r$ the indecomposable regular summands of $T$ in $\mathcal{T}_{\lambda}$. By Lemma \ref{union de conos}, we know that there are summands $E_1,...,E_m$ of $T$ such that every other summand $T_i \in \mathcal{C}(E_j)$ for some $j$ and the cones $\mathcal{C}(E_j), j=1,...,m$ are maximal respect this property. It is, $\mathcal{C}(E_j)\bigcap \mathcal{C}(E_i)= \emptyset$ if and only if $i \neq j$.

\vspace{0.2in}

If $X \in \mathcal{C}(E_j)$ for some $j$, we set $E = E_j$. Otherwise, suppose $X$ does not belong to any of this cones. We will show that can be Hom$_H(\mathcal{C}(E_{j}),X)\neq 0$ for only one of the  maximal cones $\mathcal{C}(E_{j}).$
\vspace{0.2in}

Assume that Hom$_H(\mathcal{C}(E_{j_1}),X)\neq 0$ and Hom$_H(\mathcal{C}(E_{j_2}),X)\neq 0$. Since Hom$_H(\mathcal{C}(E_{j_1}),X)\neq 0$ there exists an indecomposable module $ N_1 \in \mathcal{C}(E_{j_1}) $ in the mouth of the stable tube $\mathcal{T}_{\lambda}$ such that $X$ belongs to the ray starting in $N_1$. Analogously,  since Hom$_H(\mathcal{C}(E_{j_2}),X)\neq 0$ there exists an indecomposable module  $N_2 \in \mathcal{C}(E_{j_2}) $ in the mouth of the stable tube $\mathcal{T}_{\lambda}$ such that $X$ belongs to the ray starting in $N_2$. Then $N_1 = N_2$ because there is a unique ray in  $\mathcal{T}_{\lambda}$ containing $X$. Therefore $X \in \mathcal{C}(E_{j_1}) \bigcap \mathcal{C}(E_{j_2}) \neq \emptyset$ and by Lemma\ref{union de conos} $\mathcal{C}(E_{j_1})=\mathcal{C}(E_{j_2})$.
\vspace{0.2in}

Since $X \in \mathcal{T}(T)$ we have that Hom$_H(\mathcal{C}(E_{j}),X)\neq 0$ for some $j$. Let $E$ be the direct summand of $T$ such that Hom$_H(\mathcal{C}(E),X)\neq 0$. Let $N \in \mathcal{C}(E)$ be an indecomposable module in the mouth of the stable tube $\mathcal{T}$ such that $X$ is contained in the ray starting in $N$. We will call this ray $\mathcal{R}$. Since $X \in \mathcal{T}(T)$ there exists an epimorphism $h:T^{n} \to X$. We are assuming that $X$ does not belong to the cone $\mathcal{C}(E)$. Using the lifting property  of the Auslander-Reiten sequences the morphism $h$ must factorize by a module $C$ in the mouth  of the stable tube $\mathcal{T}_{\lambda}$, by a module $M_{1}$ in the edge of the cone $\mathcal{C}(E)$, such that there exists an epimorphism $E \to M_{1}$, or a combination of both of this cases. More over, using the mesh relations of the Auslanter-Reiten quiver, specifically  the mesh of $\mathcal{T}_{\lambda}$ we have:

\vspace{0.2in}
\begin{itemize}
\item [(a)] If $h$ factors through the  mouth of $\mathcal{T}_{\lambda}$, we may assume that $C=N$ since $X$ belongs to the ray starting in $N$. By \cite[Lema 3.7]{HRS} $h$ must factor through $\mathcal{R}$, but since $X \notin \mathcal{C}(E)$, there exists  a module $M$ in $\mathcal{C}(E\to)$ such that $M \in $Gen $(E)\bigcap\mathcal{R}$.
\item [(b)] If $h$ factors through a module in the edge of $\mathcal{C}(E)$, it is, Gen$(E)\bigcap \mathcal{C}(E)$, then $h$ should factor through a module in the same ray of $X$.
\end{itemize}

\vspace{0.2in}

We conclude that in any case the morphism $h$ must factor through $M$ and therefore $h$ applies in $X$ through  a composition of morphisms in the ray $\mathcal{R}$. Which is a contradiction, since $h$ is an epimorphism and  all morphisms morphisms in the ray are monomorphism. Hence, $X \in \mathcal{C}(E)$. $\Box$

\vspace{0.2in}

Finally, we will apply the last theorem to prove that $\mathcal{G'}$ is a summand of $M'$. In fact, suppose that $G'$ is an indecomposable direct summand of $\mathcal{G'}$. then $G'=\Hom \tilde{G})$ with $G $ an indecomposable $H$-module  such that $G \in \mathcal{T}(T)$. Therefore we have only two possibilities for  $G$, it is a regular or it is a preinjective module.

\vspace{0.1in}

If $G$ is preinjective, is easy to see from the choose of the complete slice $\Sigma$ that $G \in \mathcal{T}(\Sigma)$ and therefore $G'$ is a  direct summand of $M'_1$.

\vspace{0.1in}

If $G$ is regular, by Theorem \ref{prop torsion}, there exists a direct summand of $T$, $E$  such that $G \in \mathcal{C}(E)$. Therefore $G \in  \mathcal{C}(E)$. Then it follows from the definition of $M'_2$ that $G'$ results a direct summand of $M'$.

\vspace{0.2in}
Note that in the proof of theorem \ref{prop torsion} we have proved the following lemma:

\begin{lem}\label{borde del tubo}
Let $M \in \mathcal{C}(E)$ and $N \notin \mathcal{C}(E)$. Every non zero morphism from $M \to N $ factors trough  $\mathcal{C}(E \to)$.
\end{lem}

\vspace{0.2in}

\section{Weak representation dimension of tame cluster tilted algebras}

\vspace{0.2in}

For each $B$-module $X' \notin $ add$(M')$, we will call $M'_{X}$ to the $B$-module $P'_{X/tX} \oplus (tX)'$. In the previous section, we have constructed an epimorphism $f':M'_X \to X'$. We will begin this section studying the properties of this morphism. In particular, we want to show that the following result holds:

\vspace{0.1in}
\begin{teo}\label{fact}
The epimorphism $f':M'_X \rightarrow X'$ constructed in Proposition \ref{resolucion} is an  add$(M')$-approximation de $X'$.
\end{teo}

\vspace{0.1in}

We will concentrate for the rest of this section in proving this result. In order to do this, we will study the behavior of the morphism $f'$ in relation of the functor
Hom$_B(M'',$ \, $)$, where $M''$ is in add$(M')$.

\vspace{0.1in}

We begin by showing the following lemma.

\begin{lem}\label{lem M_X}
The morphism $f':M'_X \to X'$ is an add$(H' \oplus \mathcal{T}')$-approximation for every $X' \notin$ add$(M')$.
\end{lem}

\noindent \textbf{proof:}

We are going to split this proof in the following two cases:

\begin{itemize}
\item[a)] any morphism of $B$-modules from add$(H')$ to $X'$ factors through $f'$.
\item[b)] any morphism of $B$-modules from add$(\mathcal{T}')$ to $X'$ factors through $f'$.
\end{itemize}

\vspace{0.1in}
We start with the first case. Let $P'$ be in add$(H')$ and let $ 0 \neq h' : P' \rightarrow X'$. Without lose of generality, we may assume that $P'=$Hom$_\mathcal{C}(\tilde{T},\tilde{P})$ with $P$ an $H$-module projective indecomposable. Then, there exists $\tilde{h} \in $ Hom$_\mathcal{C}(\tilde{P},\tilde{X})$ such that $h'=$Hom$_\mathcal{C}(\tilde{T},\tilde{h})$.

\vspace{0.1in}

In fact, we can identify $\tilde{h}$ with a morphism $t$ of $H$-modules. Since $P$ is projective and $X$ in mod $H$ we have that Hom$_{\mathcal{D}}(P,\tau^{-1}X[1])=0$, therefore Hom$_\mathcal{C}(\tilde{P},\tilde{X})=$ Hom$_\mathcal{D}(P,X)=$ Hom$_H(P,X)$. Let $t:P \to X$ such that $\tilde{t}=\tilde{h}$.

\vspace{0.2in}

Consider  the canonical exact sequence in mod $H$ and the composition of the morphisms $p$ and $t$:

\vspace{0.1in}

$$\xymatrix{ & & P\ar[d]^{t} & & \\
0\ar[r] & tX\ar[r]^{i} & X \ar[r]^{p} & (X/tX)\ar[r] & 0}$$

\vspace{0.1in}


\vspace{0.2in}

Since $P$ is projective, the morphism $pt: P \to X/tX$ must factor  through the minimal projective cover of $X/tX$, it is, factors through $\pi:P_{X/tX} \to X/tX$,  the projective cover of $X/tX$, which implies that there exists $j:P \to P_{X/tX}$  such that $pt=\pi j$.

\vspace{0.2in}

Recall that by the construction of $g$, $\pi =pg$. Hence, we have that $pt = pgj$. Therefore $p(t-gj)=0$ and by the kernel universal property the morphism $t-gj$ must factors through Ker$p = tX$. Hence,  there exists $k:P \to tX$ such that $t-gj = ik$. Thus, we have proof that $t= gj + ik = f \left(\begin{array}{c}k\\ j\end{array}\right)$, it is, $t$ factors through $M_X$ in mod $H$.

\vspace{0.1in}

We define $l:= \left(\begin{array}{c}k\\ j\end{array}\right)$ and let $l'=$  Hom$_{C}(\tilde{T}, \tilde{l})$. Then is easy to see that $h'=f'l'$. Note that $l'$ must be a non zero morphism in mod $B$, since otherwise $h'$ it would be null. This completes the proof of lemma corresponding to the case (a).

Lets focus on case (b). Let $G$ be in add$(\mathcal{T}')$, we want to show that $h': G' \rightarrow X'$ it factors trough $f':M_X' \rightarrow X'$.

\vspace{0.2in}

Suppose that we have  $h':G' \rightarrow X'$. Then exists a morphisms
\newline$\tilde{h} \in \mbox{Hom}_\mathcal{C}(\tilde{G}, \tilde{X})$ such that  $h' = \Hom \tilde{h})$. We know that $$\mbox{Hom}_\mathcal{C}(\tilde{G}, \tilde{X}) = \mbox{Hom}_\mathcal{D}^{b}(H)(G, X) \oplus \mbox{Hom}_\mathcal{D}^{b}(H)(G, \tau ^{-1}X[1]).$$ Let $\tilde{h}=(h_1, h_2):\tilde{G} \to \tilde{X}$ with $h_1:G \to X$ y $h_2:G \to \tau^{-1}X[1]$. It suffices to show that this morphism factors through $\tilde{f}=( f, 0)$ where $f:M_X \to X$ in mod$H$. It is, there exists $\tilde{g}=(g_1, g_2)$ such that $\tilde{g}\tilde{f}=\tilde{h}$.

\vspace{0.2in}

We will show that there exists a morphism $\tilde{h}= (h_1, h_2)$ such that the following diagram results commutative:

$$ \xymatrix{ M_X \oplus FM_X \ar[rr]^{\left(\begin{array}{cc}f & 0\\ 0 &Ff \end{array}\right)}& & X \oplus FX \\
G\ar[u]^{(g_1, g_2)}\ar[rru]_{(h_1,h_2)}&&}$$

\vspace{0.2in}

Then we will have that $\tilde{g}\tilde{f}=(g_1, g_2)\left(\begin{array}{cc}f & 0\\ 0 &Ff \end{array}\right)$ $= (g_1f, Ffg_2)$, so it suffices to prove that any morphism $h$ en $\mathcal{D}^{b}(H)$, such that $h: G \rightarrow X $ it factors trough $f$ in $\mathcal{D}^{b}(H)$ and any morphism $h: G \rightarrow \tau^{-1}X[1]$  it factors through $Ff$ in $\mathcal{D}^{b}(H).$

\vspace{0.2in}

Suppose $h: G \rightarrow X $. Since $G \in \mathcal{T}(T)$, $\mbox{Hom}_{H}(G, X) = \mbox{Hom}_{H}(G, tX)$ then we have that $h$ factors through $tX$ trivially. Therefore, $h'$ factors trivially over $(tX)'$ which is a summand of $M'_X$. This concludes the proof for this case.

\vspace{0.2in}

Now, assume that we have a morphism in $\mathcal{D}^{b}(H)$, $h: G\rightarrow FX$. We will show that this case can be reduced to case (a).

\vspace{0.2in}

We start with a projective resolution of $X$ in mod $H$

$$0 \to K_X \to P_X \stackrel{\pi_X}\to X \to 0,$$

\vspace{0.1in}
\noindent This exact sequence induces a triangle in $\mathcal{D}^{b}(H)$:

$$ K_X \to P_X \stackrel{\pi_X}\to X \to K_X[1] $$

\vspace{0.2in}

Now, applying the functor $F$ to this triangle,  we obtain another triangle:

\vspace{0.2in}

$$ FK_X \to FP_X \stackrel{F\pi_X}\to FX \to FK_X[1] $$

\vspace{0.2in}

Finally, we apply the functor Hom$_{\mathcal{D}}(G, \, \, \,)$ to the last triangle to obtain the long exact sequence:

$$\cdots \rightarrow \mbox{Hom}_{\mathcal{D}}(G, FP_X) \stackrel{(F\pi_X)^*}\rightarrow \mbox{Hom}_{\mathcal{D}}(G, FX) \to \mbox{Hom}_{\mathcal{D}}(G, FK_X[1]) \to \cdots $$

\vspace{0.2in}
\noindent where the morphism $(F\pi_X)^*$ denotes the composition by $(F\pi_X)$.

\vspace{0.2in}
Note that $\mbox{Hom}_{\mathcal{D}}(G, FK_X[1])=\mbox{Ext}^2_H(G,\tau^{-1}K_X)=0$ since $H$ is an hereditary algebra. Then, the morphism $(F\pi_X)^*$ is an epimorphism. Therefore, since
$h: G \rightarrow FX $, there exists $k:G \rightarrow FP_X$ such that $h=(F\pi_X)k$.

\vspace{0.2in}

 We have prove that $h$ factors through $FH$ and we have a morphism $Fq:FP_X \to FX$ where $P_X \in$ add$(H')$. Then proceeding similarly to the case (a) over $q:P_X \to X$ we have that $q'$ and therefore $h'$ factors through $M'_X$. $\Box$


\vspace{0.2in}

This lemma will be helpful to show that the morphism $f'$ is an add$(M')$-approximation  for every $B$-module $X'$.

\begin{prop}\label{fact 2}
Let $X' \notin $add$(M')$ and let $f':M'_X \to X'$. Then any morphism of add$(M')$ to $X'$ it factors through de $f'$.

\end{prop}

\noindent \textbf{Proof:}

We will split the proof in two cases:

\vspace{0.1in}
\begin{itemize}
\item [a)] Any morphism from add$(M'_1)$ to $X'$ factors through $f'$.
\item [b)] Any morphism from add$(M'_2)$ to $X'$ factors through $f'$.
\end{itemize}

a)
Let $Y'$ be a module indecomposable in add$(M'_1)$. Consider $H_1$ the algebra hereditary derived equivalent to $H$, given by identifying the complete slice $\Sigma$ in mod $H$ with the complete projective slice mod$H_1$. We can assume that $Y'=$ Hom$_\mathcal{C}(\tilde{T},\tilde{Y})$ with $Y$ an $H_1$-module indecomposable in $\mathcal{F}(FH)$ identifying $FH$ with the corresponding tilting module in mod $H_1$.

\vspace{0.2in}

Let $ 0 \neq h' : Y' \rightarrow X'$. Then, there exists $\tilde{h} \in $ Hom$_\mathcal{C}(\tilde{Y},\tilde{X})$ such that $h'=$ Hom$_\mathcal{C}(\tilde{T},\tilde{h})$. Since $X' \notin$ add$(M')$, we may  assume that $X$ is a module in $\mathcal{F}(\Sigma)\subset$ mod$H$. Therefore Hom$_H(Y,X)=0$ and
Hom$_\mathcal{C}(\tilde{Y},\tilde{X})= $ Hom$_\mathcal{D}(Y,\tau^{-1}X[1])$.

\vspace{0.2in}

Let  $h:Y \to \tau^{-1}X[1]$. As we see before, we know that it is possible to identify $Y$ with a preprojective $H_1$-module  with $Y \in  \mathcal{F}(FH)$.

\vspace{0.2in}

Since we suppose that $X' \notin$ add$(M')$ then we can assume that $X \in \mathcal{F}(\Sigma)$. We have two cases to analyze: $\tau^{-1}X \notin$ add$(\Sigma)$ and then $\tau^{-1}X \in \mathcal{F}(\Sigma)$ or $\tau^{-1}X \in$ add$(\Sigma)$.

\vspace{0.2in}

If $\tau^{-1}X \notin$ add$(\Sigma)$ it is also possible to identify $\tau^{-1}X[1]$ with an $H_1$-module. More over, we have that $\tau^{-1}X[1] \in \mathcal{T}(FH)$ as $H_1$-module and that $Y \in  \mathcal{F}(FH)$ therefore, since $FH$ is a complete slice in mod $H_1$ any morphism from $Y$ to $\tau^{-1}X[1]$, in particular $h$ must factor through $FH$. Thus, the morphism $h'$ in mod $B$ factors through a $B$-module in add$(H')$. Then the proof of  this cases can be deduced from case (a) of Lemma \ref{lem M_X}.

\vspace{0.2in}

It only left to see what happens with the case $\tau^{-1}X \in$ add$(\Sigma)$. Suppose that $Y \in $ add$(\Sigma)$, then  Hom$_\mathcal{D}(Y,\tau^{-1}X[1])=$ Ext$_{H}^{1}(Y,\tau^{-1}X)=0$ since $\Sigma$ is a tilting module and therefore Ext$_{H}^{1}(\Sigma,\Sigma)=0$. Then we can assume that $Y \notin $ add$(\Sigma)$. Finally, defining $H_2$ as $\tau_{\mathcal{D}}^{-1}H_1$ we proceed analogously to the previous case. This finishes the proof of case a).

\vspace{0.2in}

b)
Let $X'$ be an indecomposable $B$-module  such that $X' \notin $ add$(M')$. Let $Y'$ be an indecomposable module in add$(W'_i)$ for some $i$. We write $X'=$ Hom$_{\mathcal{C}}(\tilde{T},\tilde{X})$ and $Y'=$ Hom$_{\mathcal{C}}(\tilde{T},\tilde{Y})$.
We may assume that $Y \notin \mathcal{T}(T).$ Recall that, by the definition of $W_i$, There exists a regular summand of the tilting module $T$, which we denote by $T_{i}$ such that $Y$ must be a regular module in the cone $\mathcal{C}(T_i)$.

\vspace{0.1in}

As we did before, we have two cases to analyse:  $h:Y \to X$ or $h:Y \to \tau^{-1}X[1]$.

\vspace{0.2in}

We will show that any morphism $h: Y \to X$ factors necessarily  by $\mathcal{T}(T)$ and therefore, as we prove in the proof of case (b) of Lemma \ref{lem M_X}, by $M_X$.

\vspace{0.2in}

 Since $X' \notin$ add$(M')$ then $X$ does not belong to the cone $\mathcal{C}(T_i)$. Therefore, by corollary \ref{borde del tubo} we have that  the morphism $h$ must factor through a module  in $\mathcal{C}(T_i \to) \subset \mathcal{T}(T)$. Then, the proof of this case follows from case (b) of the previous  Lemma.

\vspace{0.2in}

All we need to analyze the case in which $h:Y \to \tau^{-1}X[1]$. Proceeding analogously to the case (b) of the previous lemma we can see that in this case $h$ must factors through $FH$ and therefore the proof of this case follows from case a) of Lemma \ref{lem M_X}. This completes the proof of the proposition. $\Box$


\vspace{0.2in}

Summarizing,we have proved the following result:

\begin{teo}\label{last teo}
Let $H$ be a tame hereditary algebra of infinite representation type. Let $T=T_{I}\oplus T_{R}$ a tilting module with $T_{I}$ a  non zero preinjective module. Let $B$ the cluster tilted algebra given by $B=$End$_{\mathcal{C}}(\tilde{T})$. Then w.rep.dim $B=3$.
\end{teo}

\noindent \textbf{Proof:}

As a direct consequence of propositions \ref{resolucion} and \ref{fact 2} we have that gl.dim End$_{B}(M')\leq 3$. Since the module $M'$ is generator and $B$ is of infinite representation type, we conclude that w.rep.dim $B=3$. $\Box$

\begin{cor}
Let $B$ a tame cluster tilted algebra. Then w.rep.dim$B = 3$.
\end{cor}

\noindent \textbf{Proof:}

Suppose that $B$ is a tame cluster tilted algebra. Then there exists an hereditary algebra $H$ of tame representation type and a tilting module $T$ such that $B=$End$_{\tilde{C}}(\tilde{T})$. Since $H$ is tame, we know that $T$ has at least a non zero non regular summand,   otherwise the algebra $B$ would be wild. Suppose that $T$ has no non zero regular summands then by theorem \cite[theorem 3.6]{GT} repdim $B=$w.rep.dim$B=3$.

Suppose $T$ has at least one regular summand, as we discuss previously there is no lose of generality in assuming that the non regular summands of $T$ are all preinjective. By Lemma \ref{3.2*} this implies that $\mathcal{T}(T)$ results finite. Thus, by theorem \ref{last teo}  w.rep.dim$B=3$. $\Box$


\vspace{0.2in}

\subsection{About representation dimension of tame cluster tilted algebras}

The module $M'$ utilized  in the previous section to compute the weak representation dimension is not necessarily a cogenerator for mod $B$. We will illustrate this fact in the following example.

\begin{ej}
Let $T'$ be the tilting module  given in  example \ref{ej2}. Consider the cluster tilted algebra given by  $B=$End$_{\mathcal{C}}(\tilde{T'})$. We will show that there exists an injective indecomposable $B$-module  which is not a summand of $M'$.

\vspace{0.2in}

First, we observe that $S_3$ is the only indecomposable regular summand of $T'$, and it lies in the mouth of the stable tube of rank $3$ in $\Gamma_H$:

\[
\begin{array}{ccccccccccccccc}
\cdots S_3 &  &  &  & \tau^{2} S_3 &  &  &  & \tau S_3 &  &  &  & S_3 &  &  \\
& \searrow &  & \nearrow &  & \searrow &  & \nearrow &  & \searrow
&  &
\nearrow &  & \searrow &  \\
&  & \bullet &  &  &  & \bullet &  &  &  & \bullet &  &  &  & \bullet \cdots \\
& \nearrow &  & \searrow &  & \nearrow &  & \searrow &  & \nearrow
&  &
\searrow &  & \nearrow &  \\
\cdots \bullet &  &  &  & \bullet &  &  &  & \bullet &  &  &  & \bullet &  &  \\
& \searrow &  & \nearrow &  & \searrow &  & \nearrow &  & \searrow
&  &
\nearrow &  & \searrow &  \\
&  & \bullet &  &  &  & \bullet &  &  &  & \bullet &  &  &  & \bullet \cdots \\

&  & \vdots &  &  &  & \vdots &  &  &  & \vdots &  &  &  &
\end{array}
\]\vspace{.2in}

Applying the Hom$_{\mathcal{C}}(T,$ \; \; \;$)$ we obtain a description of the corresponding component in $\Gamma_B$ by suppressing  the vertex associated to $\tau S_3$ in $\Gamma_H$:

\vspace{0.2in}
\[
\begin{array}{ccccccccccccccc}
\cdots (S_3)' &  &  &  & (\tau^{2} S_3)' &  &  &  &  &  &  &  & (S_3)' &  &  \\
& \searrow &  & \nearrow &  & \searrow &  &  &  &
&  &
\nearrow &  & \searrow &  \\
&  & \bullet &  &  &  & \bullet &  & \cdots &  & \bullet &  &  &  & \bullet \cdots \\
& \nearrow &  & \searrow &  & \nearrow &  & \searrow &  & \nearrow
&  &
\searrow &  & \nearrow &  \\
\cdots \bullet &  &  &  & \bullet &  &  &  & \bullet &  &  &  & \bullet &  &  \\
& \searrow &  & \nearrow &  & \searrow &  & \nearrow &  & \searrow
&  &
\nearrow &  & \searrow &  \\
&  & \bullet &  &  &  & \bullet &  &  &  & \bullet &  &  &  & \bullet \cdots \\

&  & \vdots &  &  &  & \vdots &  &  &  & \vdots &  &  &  &
\end{array}
\]\vspace{.2in}

 Note that the $B$-module $(\tau^{2}S_3)'$ is an indecomposable injective module,  and it is not transjective  in
 mod $B$. More over, the only indecomposable summand of $M_2'$, is the projective $B$-module  corresponding to vertex 3, $(S_3)'$. therefore, in this case, $M'$ is not a cogenerator for mod $B$.

\end{ej}

\vspace{0.2in}

\begin{cor}
Let $B$ a tame cluster tilted algebra.
If every non transjective injective module in $B$ lies in add$(M'_2)$ then rep.dim $B=3$.
\end{cor}

\noindent \textbf{Proof:}

Assume  that every indecomposable injective module which are not transjectives  lies in add$(M'_2)$. Recall that every injective module in $B$ is of the form $\Hom \tau^{2}\tilde{T})$. It is possible to choose $\Sigma$ in mod $H$ such that every preinjective module of the form $\tau^2T$ are in $\mathcal{T}(\Sigma)$. This means that it is possible to construct $M'_1$ with this $\Sigma$. This implies that every preinjective injective $B$-module is in add$(M'_1)$. Hence $M'$ is a cogenerator for mod $B$. Thus, in this case, follows as an immediately consequence of our last theorem that rep.dim $B=3$.$\Box$

\vspace{0.2in}

Finally, we illustrate the situation of last corollary in the following example:

\begin{ej}
Let $H$ be the hereditary algebra given by the quiver $\mathbb{\tilde{D}}_4$.
Consider the tilting $H$-module $T = I_5 \oplus I_4 \oplus I_1 \oplus I_2 \oplus S_3$. Here every indecomposable  summand of $T$ are preinjectives with the exception of $S_3$ and is immediate to see that the torsion associated to this module is finite. Let $B$ given by End$_{\mathcal{C}}(\tilde{T})$.
We know that $S_3$ is a regular module in an homogeneous stable tube in $\Gamma_H$. Then $\tau^{2}S_3=S_3$. Therefore the module $M'_2=(S_3)'$ is a projective-injective module in mod $B$. More over, choose $\Sigma=\tau^{2}DH$. Hence, by the last corollary  we conclude that rep.dim $B=3$.
\end{ej}



\begin{thebibliography}{Dillo 83}

\bibitem[ABS] {ABS} I. Assem, T. Br\"{u}stle , R. Schiffler, , Cluster-tilted algebras and slices, Journal of Algebra, vol. 319, no. 8, pág. 3464-3479, 2008.
\bibitem[AK] {AK} I. Assem, O. Kerner, Constructing torsion pairs, Journal of  Algebra, vol.185, pág.  19-41, 1996.
\bibitem[APT] {APT} I. Assem, M. I. Platzeck, S. Trepode, On the representation dimension of tilted and laura algebras, Journal of Algebra, vol.296, p\'ag. 426-439, 2006.
\bibitem[ASS] {ASS} I. Assem, A. Skowronski, D. Simpson, Elements of representation theory of associative algebras, vol 1: Techniques of representation theory, London Mathematical Society Student Texts 65, Cambridge University Press, Cambridge, 2006.
\bibitem[Au] {Au} M.Auslander, representation dimension of artin algebras, Queen Mary College Mathematics Notes, London, 1971.
\bibitem [APR]{APR} M. Auslander, M. Platzeck, I. Reiten, Coxeter functors without diagrams, Transactions of the American Mathematical Society, vol. 250, pag. 1-46, 1979.
\bibitem [AR] {AR} M. Auslander, I. Reiten, Applications of contravariantly finite subcategories, Advances in mathematics, vol. 86, Issue 1, pag. 111-152, 1991.
\bibitem [ARS]{ARS} M. Auslander, I. Reiten, S.Smal{\o}, Representation theory of artin algebras, Cambridge Studies in Advanced Mathematics, Cambridge University Press, 1997.
\bibitem [AS]{AS} M. Auslander, S. Smal{\o}, Postprojective modules over artin algebras, Journal of Algebra, vol. 66, Issue. 1, pag. 61-122, 1980.
\bibitem[BMRRT] {BMRRT} A. Buan, R. Marsh, M. Reineke, I. Reiten, G. Todorov, Tilting theory and cluster combinatorics, Advances in Mathematics, vol. 204, Issue. 2, pag. 572-618, 2006.
\bibitem[BMR1] {BMR1} A. Buan, R. Marsh, I. Reiten, Cluster-tilted algebras of finite representation type, Journal of Algebra, vol. 306, Issue 2, pag. 412-431, 2006.
\bibitem[BMR2] {BMR2} A. Buan, R. Marsh, I. Reiten, Cluster-tilted algebras, Transactions of the American Mathematical Society, vol.359, Issue 1, pag. 323-332, 2007.
\bibitem[BMR3] {BMR3} A. Buan, R. Marsh, I. Reiten, Cluster mutation via quiver representations, Commentarii Mathematici Helvetici,  vol. 83, Issue 1 , pag. 143-177, 2008.
\bibitem [BRS]{BRS} A. Buan, I. Reiten, A. Seven, Tame concealed algebras and cluster quivers of minimal infinite type, Journal of pure and applied algebra, Vol. 211, Issue. 1, pag. 71-82, 2007.
\bibitem[CP] {CP} F. Coelho, M. Platzeck, On the representation dimension of some classes of algebras, Journal of Algebra, vol. 275, pag.  615-628, 2004.
\bibitem[EHIS] {EHIS} K. Erdmann, T. Holm, O. Iyama, J. Schr\"{o}er, Radical embeddings and representation dimension,  Advances in mathematics, vol. 185, Issue. 1, pag. 159-177, 2004.
\bibitem[FZ]{FZ} S. Fomin, A. Zelevinsky, Cluster Algebras I: Foundations, Journal of the American Mathematical Society, vol. 15, Issue. 2, pag. 497-529, 2002.
\bibitem [H]{H} D. Happel, Triangulated Categories in the representation theory of Finite Dimensional Algebras, London Mathematical Society Lecture Notes 119, 1988.
\bibitem [GS]{GS} A. Garc\'{\i}a Elsener, R. Schiffler, On syzygies over 2-Calabi-Yau tilted algebras ,arXiv:1601.03988v2
\bibitem [GT]{GT} A. Gonz\'alez Chaio, S. Trepode, Representation dimension of cluster concealed algebras, Algebras and Representation Theory, Vol. 16, Issue 4, pag 1001-1015, 2013
\bibitem[HRS]{HRS} D. Happel, I.Reiten, S.Smal{\o} tilting in abelian categories and quasi-tilted algebras, Memoirs of the American Mathematical Society, vol. 120, 1996.
\bibitem [HR]{HR} D. Happel, C. Ringel, Tilted algebras, Transactions of the American Mathematical Society, vol. 274, Issue. 2, pag. 399-443,1982.
\bibitem [IT]{IT} K. Igusa, G. Todorov, On the finitistic global dimension conjecture for Artin algebras, Representations of Algebras and Related Topics, Fields Inst. Commun.,American Mathematical Society, vol. 45, pag. 201-204, 2005.
\bibitem[I]{I} O. Iyama, Finiteness of representation dimension. Proceedings of the american mathematical society, vol.131, Issue. 4,  pag. 1011-1014, 2003.
\bibitem[K]{K} B. Keller,  Triangulated orbit categories, Documenta Mathematica,  vol.10, pag. 551-581, 2005.
\bibitem[Ri1] {Ri1} C. Ringel, Tame Algebras and Integral Quadratic Forms, Lecture Notes in Mathatematics, vol. 1099, 1984.
\bibitem[Ri2] {Ri2} C. Ringel, The regular components of the Auslander-Reiten quiver of a tilted algebra, Chinese Annals of
Mathematics, vol. 9 Issue. 1, pag. 1-18, 1988.
\bibitem [R]{R} R. Rouquier, On the representation dimension of exterior algebras. Inventiones mathematicae, vol. 165, pag. 357-367, 2006.
\bibitem[SS1] {SS}  A. Skowronski, D. Simpson,  Elements of the Representation Theory of Associative Algebras 2: Tubes and Concealed Algebras of Euclidean Type, London Mathematical Society Student Texts 71, Cambridge University Press, 2007.
\bibitem[SS2]{SS2}  A. Skowronski, D. Simpson, Elements of the Representation Theory of Associative Algebras 3: Representation-Infinite Tilted Algebras, London Mathematical Society
Student Texts 72, Cambridge University Press, 2007.


\end{thebibliography}
\end{document}